\documentclass[a4paper,11pt]{article}
\usepackage{amsmath,amsthm,amssymb}
\theoremstyle{plain}
\newtheorem{theorem}{Theorem}
\newcommand{\dd}{\mathrm{d}}
\newcommand{\ee}{\mathrm{e}}
\newcommand{\ii}{\mathrm{i}}
\newcommand{\diag}{\mathop{\rm diag}}
\renewcommand{\Im}{\mathop{\rm Im}}
\newcommand{\ord}{\mathrm{O}}
\begin{document}
\title{{\Large\textbf{The mKdV equation on a finite interval}}}
\author{{\normalsize Anne \textsc{Boutet} de \textsc{Monvel}$^\ast$
         and Dmitry \textsc{Shepelsky}$^\dagger$}\\[1mm]
{\scriptsize $^{\ast}$ Institut de Math\'ematiques de Jussieu, case 7012, 
     Universit\'e Paris 7,}\\[-2mm] 
{\scriptsize 2 place Jussieu, 75251 Paris, France}\\
{\scriptsize 
$^{\dagger}$ Mathematical Division, Institute for Low Temperature Physics,}\\[-2mm]
{\scriptsize 47 Lenin Avenue, 61103 Kharkiv, Ukraine}}
\date{}
\maketitle
\begin{abstract}
We analyse an initial-boundary value problem for the mKdV equation
on a finite interval 
by expressing 
the solution in terms of the solution of an associated matrix Riemann-Hilbert problem
in the complex $k$-plane. 
This Riemann-Hilbert problem has explicit $(x,t)$-dependence and it
involves certain functions of $k$ referred to as ``spectral functions''.
Some of these functions are 
defined in terms
of the initial condition $q(x,0)=q_0(x)$, while the
remaining spectral functions are defined in terms of two sets of boundary values.
We show that the spectral functions satisfy an algebraic ``global relation''
that characterize the boundary values in spectral terms.
\end{abstract}
\section{Introduction}

The general method 
for solving initial-boundary value problems
for two-dimensional linear and integrable nonlinear PDEs
announced in \cite{F97} and developed further in 
\cite{F97}\textendash\cite{F01} 
is based on the simultaneous spectral
analysis of the two eigenvalue equations of the associated Lax pair.
It expresses  the solution in terms of the solution of a matrix
Riemann-Hilbert (RH) problem  formulated in the complex plane of the 
spectral parameter. The spectral functions determining the RH problem are
expressed in terms of the initial and boundary values of the solution. The
fact that these values are in general related can be expressed in a
simple way in terms of a global relation satisfied by the corresponding
spectral functions.

The rigorous implementation of the method to the 
modified
Korteweg--de Vries (mKdV) equation  on the half-line is presented in 
\cite{BFS}. In
the present Note, this methodology is applied  to 
the mKdV equation on a finite interval. The similar problem for the nonlinear 
Schr\"odinger
equation is studied in \cite{FI}.

The modified Korteweg\textendash de Vries equation 
\begin{equation}\label{mkdv}
q_t - q_{xxx} + 6\lambda q^2 q_x = 0, \quad \lambda =\pm 1
\end{equation}
admits the Lax pair formulation
\begin{equation}\label{lax}
\mu_x -\ii k\hat\sigma_3\mu = Q(x,t)\mu, \qquad
\mu_t + 4\ii k^3\hat\sigma_3\mu = \tilde Q(x,t,k)\mu,
\end{equation}
where $\sigma_3=\diag\{1,-1\}$, 
$\hat\sigma_3A:= \sigma_3A - A\sigma_3$,
$\ee^{\hat\sigma_3}A=\ee^{\sigma_3}A\ee^{-\sigma_3}$,
\[
Q(x,t) = 
\begin{pmatrix}
0 & q(x,t) \\ 
\lambda q(x,t) & 0
\end{pmatrix}, \quad \tilde Q(x,t,k) = -4k^2Q-2\ii k(Q^2+Q_x)\sigma_3
-2Q^3+Q_{xx}.
\]
We study  the initial-boundary value problem for
the mKdV equation in the domain $\{0<x<L,\, 0<t<T\}$, $L<\infty$,
$T\leq\infty$ using the following steps.\\[1mm]
\textbullet\
Assuming that the solution $q(x,t)$ of the mKdV equation exists, express
it via the solution of a matrix Riemann-Hilbert problem. For this purpose:
\begin{enumerate}
\item 
Define proper solutions of (\ref{lax})  sectionally analytic  and bounded in
 $ k\in\overline{\mathbb{C}}=\mathbb{C}\cup\{\infty\}$.
\item 
Define spectral functions $s(k)$, $S(k)$, and $S_1(k)$ such that:
\begin{itemize}
\item 
They determine a Riemann-Hilbert problem.
\item  
$s(k)$ is determined by the initial conditions $q(x,0)=q_0(x)$, $0<x<L$.
\item  
$S(k)$ is determined by the boundary values $q(0,t)=g_0(t)$, 
$q_x(0,t)=g_1(t)$, $q_{xx}(0,t)=g_2(t)$, $0<t<T$.
\item  
$S_1(k)$ is determined by the boundary values $q(L,t)=f_0(t)$, 
$q_x(L,t)=f_1(t)$, $q_{xx}(L,t)=f_2(t)$, $0<t<T$.
\item  
They satisfied an algebraic ``global relation'', expressing the fact that 
$q_0(x)$, $\{g_j(t)\}_{j=0}^2$, $\{f_j(t)\}_{j=0}^2$ being the
initial and boundary conditions for the mKdV equation, 
cannot be chosen arbitrarily.
\end{itemize}
\end{enumerate}
\textbullet\
Given $s(k)$ and assuming that $\{g_j(t)\}_{j=0}^2$ and $\{f_j(t)\}_{j=0}^2$
are such that
the associated $S(k)$ and $S_1(k)$ together with $s(k)$ satisfy the global
relation, prove that the  solution of the Riemann-Hilbert problem constructed
from $s(k)$, $S(k)$, and $S_1(k)$ generates the solution of the
initial-boundary value problem for the mKdV equation with  initial data
$q(x,0)=q_0(x)$ and  boundary values $q(0,t)=g_0(t)$,  $q_x(0,t)=g_1(t)$, 
$q_{xx}(0,t)=g_2(t)$,  $q(L,t)=f_0(t)$,  $q_x(L,t)=f_1(t)$, 
$q_{xx}(L,t)=f_2(t)$.

\section{Eigenfunctions and 
                    spectral functions}

Assume that there exists a real-valued function $q(x,t)$ with sufficient
smoothness and decay satisfying (\ref{mkdv}) in $\{0<x<L, 0<t<T\}$, $T\leq\infty$. 
Define the
\emph{eigenfunctions} $\mu_n(x,t,k)$, $n=1,2,3,4$ as matrix-valued solutions
of the integral equations
\begin{equation}
\mu_n(x,t,k) 
= I + \int_{(x_n,t_n)}^{(x,t)}\ee^{\ii(k(x-y)-4k^3(t-\tau))\hat\sigma_3}
(Q\mu_n\dd y +\tilde Q\mu_n\dd\tau),
\end{equation}
where $(x_1,t_1) = (0,T)$, $(x_2,t_2)=(0,0)$, $(x_3,t_3)=(L,0)$, 
$(x_4,t_4)=(L,T)$,
and the paths of integration are chosen to be parallel 
to the $x$ and $t$ axes:
\begin{align*}                               
\mu_1(x,t,k) 
&= 
I+\int_0^x\ee^{\ii k(x-y)\hat\sigma_3}(Q\mu_1)(y,t,k)\dd y -
\ee^{\ii kx\hat\sigma_3}
\int_t^T\ee^{-4\ii k^3(t-\tau)
\hat\sigma_3}(\tilde Q\mu_1)(0,\tau,k)\dd\tau,\\
\mu_4(x,t,k) 
&= 
I-\int_x^L\ee^{\ii k(x-y)\hat\sigma_3}(Q\mu_4)(y,t,k)\dd y-
\ee^{\ii k(x-L)\hat\sigma_3}
\int_t^T\ee^{-4\ii k^3(t-\tau)\hat\sigma_3}(\tilde Q\mu_4)(L,\tau,k)\dd\tau.
\end{align*}                               
Equations for $\mu_2$ and $\mu_3$ are similar to those for $\mu_1$ and
$\mu_4$, respectively, with the integral term $\int_0^t$ instead
of $-\int_t^T$.  
The columns of 
$\mu_n=(\begin{matrix}\mu_n^{(1)}& \mu_n^{(2)}\end{matrix})$ are
analytic and  bounded in domains  
separated by
the three lines
$\{k\in\mathbb C\mid\Im k^3=0\}$,
see Figure 1:
\[
\mu_1^{(1)}, \mu_3^{(2)}\text{ in IV}\cup\text{VI};\quad
\mu_1^{(2)}, \mu_3^{(1)}\text{ in I}\cup\text{III};\quad
\mu_2^{(1)}, \mu_4^{(2)}\text{ in V};\quad
\mu_2^{(2)}, \mu_4^{(1)}\text{ in II}.
\]

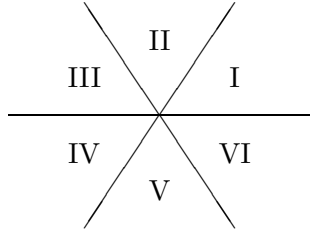
\begin{figure}[h]
\setlength{\unitlength}{5mm}
\hspace{4cm}
\begin{picture}(12,6)
\put(3,3){\line(1,0){8}}
\put(7,3){\line(2,3){2}}
\put(7,3){\line(-2,3){2}}
\put(7,3){\line(2,-3){2}}
\put(7,3){\line(-2,-3){2}}
\put(9,4){\makebox(0,0){I}}
\put(5,4){\makebox(0,0){III}}
\put(7,5){\makebox(0,0){II}}
\put(9,2){\makebox(0,0){VI}}
\put(5,2){\makebox(0,0){IV}}
\put(7,1){\makebox(0,0){V}}
\end{picture}
\caption{Domains of boundedness of eigenfunctions}
\end{figure}

\noindent
Thus, in each domain I,\dots,VI, one has a bounded $2\times 2$
matrix-valued  eigenfunction,
consisting of the appropriate vectors $\mu_n^{(l)}$.
The eigenfunctions $\mu_j$  are  related by
\begin{align} 
\mu_3(x,t,k) 
&= \mu_2(x,t,k)\ee^{\ii(kx-4k^3t)\hat\sigma_3}s(k),\label{rel-a}\\
\mu_1(x,t,k) 
&= \mu_2(x,t,k)\ee^{\ii(kx-4k^3t)\hat\sigma_3}S(k),\label{rel-b}\\
\mu_4(x,t,k) 
&=
\mu_3(x,t,k)\ee^{\ii(kx-4k^3t)\hat\sigma_3}
\ee^{-\ii kL\hat\sigma_3}S_1(k)\label{rel-c},
\end{align}
where the \emph{spectral (matrix-valued) functions} are defined as follows:
\begin{align*}
&s(k) \equiv 
\begin{pmatrix}
\overline{a(\bar k)} & b(k) \\ 
\lambda\overline{b(\bar k)} & a(k)
 \end{pmatrix}:= \mu_3(0,0,k),\\
\label{spec-2}
&S(k)\equiv 
\begin{pmatrix}
\overline{A(\bar k)} & B(k) \\ \lambda\overline{B(\bar k)} & A(k)
 \end{pmatrix}:= \mu_1(0,0,k),\\
&S_1(k)\equiv 
\begin{pmatrix}
\overline{A_1(\bar k)} & B_1(k) \\ 
\lambda\overline{B_1(\bar k)} & A_1(k)
\end{pmatrix}:= \mu_4(L,0,k).
\end{align*}
The direct and inverse spectral maps 
\begin{align*}
&\{q_0(x)\}\leftrightarrow\{a(k),b(k)\},\\
&\{g_0(t),g_1(t),g_2(t)\}\leftrightarrow\{A(k),B(k)\},\\ 
&\{f_0(t),f_1(t),f_2(t)\}\leftrightarrow\{A_1(k),B_1(k)\}
\end{align*}
are well-defined \cite{BFS}. They correspond to the separate spectral maps
for the $x$-problem ($t=0$) and $t$-problems ($x=0$ and $x=L$)
from the Lax pair (\ref{lax}).

\section{Global relation}

Evaluating equations \eqref{rel-a} and \eqref{rel-c} at $x=0$, $t=T$ and 
writing $\mu_3(0,0,k)$, $\mu_2(0,T,k)$, and $\mu_4(L,0,k)$  
  in terms of $s(k)$, $S(k)$, and 
$S_1(k)$, respectively, we obtain
\begin{equation}\label{gr1}
S^{-1}(k)s(k)\left[\ee^{-\ii kL\hat\sigma_3}
S_1(k)\right] 
=I-\ee^{4\ii k^3T\hat\sigma_3}\int_0^L\ee^{-\ii ky\hat\sigma_3}(Q\mu_4)(y,T,k)\dd y.
\end{equation}
\textbullet\
For $T<\infty$, the $(1,2)$ coefficient of (\ref{gr1}) is
($k\in\mathbb{C}$)
\begin{equation}\label{GRT}
\ee^{-2\ii kL}\left(\overline{a(\bar k)}A(k) -\lambda \overline{b(\bar
k)}B(k)\right)B_1(k)  -\left(a(k)B(k) - b(k)A(k)\right)A_1(k)=\ee^{8\ii k^3T} c(k), 
\end{equation} 
 where
$
c(k)= \int_0^L\ee^{-2\ii ky}(Q\mu_4)_{12}(y,T,k)\dd y
$
is an entire function which is $\ord\left((1+\ee^{-2\ii kL})/k\right)$ 
as $k\to\infty$.\\[1mm]
\textbullet\
For $T=\infty$, the $(1,2)$ coefficient of (\ref{gr1}) becomes
\begin{equation}\label{GRinf}
\ee^{-2\ii kL}
\left(\overline{a(\bar k)}A(k) -\lambda \overline{b(\bar k)}B(k)\right)B_1(k) -
\left(a(k)B(k)-b(k)A(k)\right)A_1(k)= 0,
\end{equation}
which is valid for $k\in\text{I}\cup\text{III}\cup\text{V}$.

Equation (\ref{GRT}) for $T<\infty$, or (\ref{GRinf}) for
$T=\infty$, is an algebraic relation between  the spectral 
functions. We call it ``global relation'', because it express, in
spectral terms, the relations between the initial and boundary values of a 
solution of
the mKdV equation. 
The global relation can be used to characterize the unknown boundary
values in a well-posed boundary value problems, say,  $g_2(t)$, $f_1(t)$, and $f_2(t)$
in terms of the boundary conditions 
$\{q_0(x), g_0(t), g_1(t), f_0(t)\}$.

\section{The Riemann-Hilbert problem}

Define a sectionally holomorphic, matrix-valued function $M(x,t,k)$:
\begin{equation}\label{M}
M=
\begin{cases}
\begin{pmatrix}
\mu_3^{(1)}& \frac{\mu_1^{(2)}}{\overline{d(\bar
k)}}
\end{pmatrix}, 
& k\in\text{I}\cup\text{III},\ |k|>R\\[3mm]
\begin{pmatrix}
\frac{\mu_4^{(1)}\overline{a(\bar k)}}{\overline{d_1(\bar k)}},
\frac{\mu_2^{(2)}}{\overline{a(\bar k)}}
\end{pmatrix}, & k\in\text{II},\ |k|>R\\[3mm] 
\begin{pmatrix}
\frac{\mu_1^{(1)}}{d(k)}, \mu_3^{(2)}
\end{pmatrix}, 
& k\in\text{IV}\cup\text{VI},\ |k|>R\\[3mm]
\begin{pmatrix}
\frac{\mu_2^{(1)}}{a(k)}, \frac{\mu_4^{(2)}a(k)}{d_1(k)}
\end{pmatrix}, 
  & k\in\text{V},\ |k|>R\\
\;\mu_2&|k|<R,
\end{cases}
\end{equation}
where 
$d(k) = a(k)\overline{A(\bar k)} -\lambda b(k)\overline{B(\bar k)}$, 
$d_1(k) = a(k)A_1(k) +\lambda\ee^{-2\ii kL}\overline{b(\bar k)}B_1(k)$,
and  $R$ is large enough
so that all possible zeros of $a(k)$, $d(k)$, and $d_1(k)$ in $\Im k\leq 0$ 
are in the disk $|k|<R$.  

Denote by $\Sigma$ the contour $\{k\mid\Im k^3=0\}\cup\{k\mid|k|=R\}$ 
(Figure 2). 

\begin{figure}[h]
\setlength{\unitlength}{5mm}
\hspace{4cm}
\begin{picture}(12,5)
\put(3,3){\line(1,0){8}}
\put(7,3){\line(2,3){2}}
\put(7,3){\line(-2,3){2}}
\put(7,3){\line(2,-3){2}}
\put(7,3){\line(-2,-3){2}}
\put(7,3){\circle{8}}
\put(9,4){\makebox(0,0){$+$}}
\put(5,4){\makebox(0,0){$+$}}
\put(7,5){\makebox(0,0){$-$}}
\put(9,2){\makebox(0,0){$-$}}
\put(5,2){\makebox(0,0){$-$}}
\put(7,1){\makebox(0,0){$+$}}
\put(7.8,3.5){\makebox(0,0){$-$}}
\put(6.1,3.5){\makebox(0,0){$-$}}
\put(7,3.8){\makebox(0,0){$+$}}
\put(7.8,2.6){\makebox(0,0){$+$}}
\put(6.1,2.6){\makebox(0,0){$+$}}
\put(7,2.2){\makebox(0,0){$-$}}
\end{picture}
\caption{Contour $\Sigma$ and domains $\Omega_\pm$}
\end{figure}
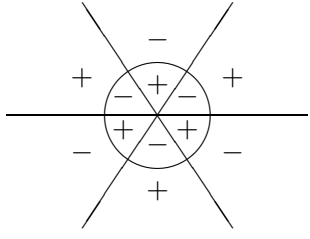
\noindent
Then the limit  values $M_\pm(x,t,k)$ (as $k$ approaches
$\Sigma$  from $\Omega_\pm$) of $M(x,t,k)$ are related on $\Sigma$ by a 
jump matrix:
\begin{equation}\label{RHglob}
M_-(x,t,k) = M_+(x,t,k) 
\ee^{(\ii kx-4\ii k^3t)\sigma_3}J_0(k)\ee^{-(\ii kx-4\ii k^3t)\sigma_3}, 
\qquad k\in\Sigma,
\end{equation}
where 
\begin{align}\label{Jglob1}
J_0(k)
&=
\begin{cases}
  \begin{pmatrix}
			1 & -\lambda\overline{\Gamma(\bar k)} \\
   0 & 1
		 \end{pmatrix}
  \begin{pmatrix}
		1 & 0 \\ \lambda\overline{\Gamma_1(\bar k)} & 1
		\end{pmatrix}, & \arg k = \frac{\pi}{3},\frac{2\pi}{3};\ |k|>R\\[3mm]
  \begin{pmatrix}
		1 & -\lambda\overline{\Gamma(\bar k)} \\  0 & 1
  \end{pmatrix}
  \begin{pmatrix}
         1-\lambda|\gamma(k)|^2 & \gamma(k) \\
                  -\lambda \bar\gamma(k) & 1
  \end{pmatrix}
  \begin{pmatrix}
    1 & 0  \\ 
    \Gamma(k) & 1
  \end{pmatrix}, & \arg k = 0,\pi;\ |k|>R \\[4mm] 
  \begin{pmatrix}
		\frac{A(k)}{\overline{d(\bar k)}} & -\dfrac{B(k)}{\overline{d(\bar k)}} \\ 
			-\lambda\overline{b(\bar k)} & \overline{a(\bar k)} 
		\end{pmatrix}, & 
					\arg k\in\Bigl(0,\frac{\pi}{3}\Bigr)\cup\Bigl(\frac{2\pi}{3},\pi\Bigr);\ |k|=R \\ 
  \begin{pmatrix}
			\overline{a(\bar k)} & 0 \\  \lambda\overline{\Gamma_2(\bar k)} &
\dfrac{1}{\overline{a(\bar k)}} 
  \end{pmatrix}, & \arg k\in\Bigl(\frac{\pi}{3},\frac{2\pi}{3}\Bigr);\ |k|=R 
\end{cases}\\
&\quad\;\text{ for }k\in\Sigma,\ \Im k\geq 0,\notag\\
J_0(k)&=\diag\{-1,\lambda\}J_0^*(\bar k)\diag\{-1,\lambda\}
\text{ for }k\in \Sigma,\ \Im k<0,\notag\\
J_0(k)&=I
\text{ for }k\in \Sigma,\ |k|<R.\notag
\end{align}
Here
\begin{align}\label{gammas}
&\gamma(k) = \dfrac{b(k)}{\bar a(k)},\notag \\
&\Gamma(k) = \lambda \dfrac{\overline{B(\bar k)}/\overline{A(\bar k)}}
			{a(k)\left(a(k)-\lambda b(k)(\overline{B(\bar k)}/\overline{A(\bar k)}\right)},
   \notag\\
&\Gamma_1(k) = \dfrac{\ee^{-2\ii kL}a(k)(B_1(k)/A_1(k))}{a(k)
+\lambda\ee^{-2\ii kL} \overline{b(\bar k)}
		(B_1(k)/A_1(k))},\notag \\
&\Gamma_2(k)= 
a(k)\dfrac{\ee^{-2\ii kL}\overline{a(\bar k)}(B_1(k)/A_1(k))+b(k)}{a(k)
+\lambda\ee^{-2\ii kL}\overline{b(\bar k)}(B_1(k)/A_1(k))}.
\end{align}
Therefore, the jump data in (\ref{RHglob}) are determined by
$a(k)$ and $b(k)$ for $k\in\mathbb{C}$, $|k|\geq R$ and by
$B(k)/{A(k)}$ and $B_1(k)/{A_1(k)}$ for 
$k\in\text{I}\cup\text{III}\cup\text{V}$, $|k|\geq R$.

\begin{theorem}\label{thmain}
Let $q_0(x)\in S(\mathbb{R}^+)$. Suppose that the sets of functions
$\{g_j(t)\}_{j=0}^2$  and $\{f_j(t)\}_{j=0}^2$ are such that the 
associated spectral
functions
$s(k)$, $S(k)$, and  $S_1(k)$ satisfy 
the global relation $(\ref{GRT})$ for $T<\infty$, or $(\ref{GRinf})$ for
$T=\infty$, where $c(k)$ is an entire function such that $c(k)=
\ord\left((1+\ee^{-2\ii kL})/{k}
\right)$ as $|k|\to\infty$. 

Let $M(x,t,k)$ be a solution of the following $2\times 2$ matrix 
RH problem:
\begin{itemize}
\item[]
\textbullet\
$M$ is sectionally holomorphic in $k\in\mathbb{C}\setminus\Sigma$.
\item[]
\textbullet\
At $k\in\Sigma$, $M$ satisfies the jump conditions $(\ref{RHglob})$, 
where  $J_0$ is defined in terms of the spectral functions
$a,b,A,B,A_1$, and $B_1$ by eqs.\ $(\ref{Jglob1})$, $(\ref{gammas})$.
\item[]
\textbullet\
$M(x,t,k) = I +\ord\left(1/{k}\right)$ as $k\to\infty$.
\end{itemize}
Then:
\begin{enumerate}
\item[\emph{(i)}] 
$M(x,t,k)$ exists and is unique;
\item[\emph{(ii)}]  
$q(x,t):= -2\ii\lim_{k\to\infty}\left(kM(x,t,k)\right)_{12}$
satisfies the mKdV equation $(\ref{mkdv})$;
\item[\emph{(iii)}] 
$q(x,t)$ satisfies the initial condition $q(x,0)=q_0(x)$ and
boundary conditions
$q(0,t)=g_0(t)$, $q_x(0,t)=g_1(t)$, $q_{xx}(0,t)=g_2(t)$, and 
$q(L,t)=f_0(t)$,
$q_x(L,t)=f_1(t)$, $q_{xx}(L,t)=f_2(t)$.
\end{enumerate}
\end{theorem}

\begin{proof}[Sketch of proof]
The  unique solvability of the RH problem is a consequence of a 
``vanishing lemma'' for
the associated RH problem with vanishing condition at infinity $M=\ord(1/k)$,
$k\to\infty$. 

The proof that the function $q(x,t)$ thus constructed solves the mKdV equation is 
straightforward and
follows the proof in the case of the whole line problem.

The proof that $q$ satisfies the initial  condition $q(x,0)=q_0(x)$
follows from the fact that it is
possible
 to map the RH problem for $M(x,0,k)$ to that for a sectionally holomorphic
function $M^{(x)}(x,k)$ corresponding to the spectral problem for
the $x$-part of the Lax pair (\ref{lax}):
$M^{(x)}(x,k) = M(x,0,k)P^{(x)}(x,k)$ 
where 
$P^{(x)}$ is sectionally holomorphic  and 
$P^{(x)}= I + P^{(x)}_{\text{off}}$, with $P^{(x)}_{\text{off}}(x,k)$  
off-diagonal 
and exponentially decaying as $k\to\infty$ for $\Im k\neq 0$.

The proof that $q$ satisfies the boundary  conditions is, in turn, 
based on the consideration of the maps $M(0,t,k)\mapsto M^{(t)}(t,k)$
and  $M(L,t,k) \mapsto M^{(t)}_1(t,k)$,
where $ M^{(t)}(t,k)$ and  $M^{(t)}_1(t,k)$ correspond
to the spectral problems for
the $t$-equation in the Lax pair (\ref{lax}) at $x=0$ and $x=L$:
$M^{(t)}(t,k)= M(0,t,k)P^{(t)}(t,k)$, $M^{(t)}_1(t,k) = M(L,t,k)P^{(t)}_1(t,k)$. 
In this case, it is the global relation (\ref{GRT}), or (\ref{GRinf}), that
guarantees that  
$P^{(t)} = P^{(t)}_{\text{diag}} + P^{(t)}_{\text{off}}$, where
$P^{(t)}_{\text{diag}}$ is  diagonal, $P^{(t)}_{\text{diag}} = I
+\ord\left(1/{k}\right)$,  and
$P^{(t)}_{\text{off}}(t,k)$  is off-diagonal and exponentially decaying 
as $k\to\infty$, and similarly for $P^{(t)}_1(t,k)$.
\end{proof}


\end{document}